\documentclass[12pt]{amsart}

\usepackage{xypic}
\usepackage{amsfonts}

\newcommand{\bintree}{\xymatrix@!R@!C@=2ex@M=0.3ex}
\newcommand{\LEFT}{\ar@{-}[dl]}
\newcommand{\RIGHT}{\ar@{-}[dr]}
\newcommand{\DOWN}{\ar@{-}[d]}
\newcommand{\LDR}{\LEFT\DOWN\RIGHT}

\newcommand{\LD}{\LEFT\DOWN}
\newcommand{\DR}{\DOWN\RIGHT}

\title[Pattern Avoidance with Additional
Restrictions]{Generalized Pattern Avoidance with Additional
Restrictions}
\author{Sergey Kitaev}
\email{kitaev@math.chalmers.se} 
\address{Matematik, Chalmers tekniska h\"ogskola och G\"oteborgs
universitet, S-412~96  G\"oteborg, Sweden}

\newtheorem{prop}{Proposition}
\newtheorem{lemma}[prop]{Lemma}

\newtheorem{thm}[prop]{Theorem}
\theoremstyle{definition}

\def\newop#1{\expandafter\def\csname #1\endcsname{\mathop{\rm
#1}\nolimits}}
\newop{erf}

\begin{document}

\begin{abstract}
Babson and Steingr\'{\i}msson introduced generalized permutation
patterns that allow the requirement that two adjacent letters in a
pattern must be adjacent in the permutation. We consider
$n$-permutations that avoid the generalized pattern $1-32$ and whose $k$
rightmost letters form an increasing subword. The number of such
permutations is a linear combination of Bell numbers. We find a
bijection between these permutations and all partitions of an
$(n-1)$-element set with one subset marked that satisfy certain
additional conditions. Also we find the e.g.f. for the number of
permutations that avoid a generalized 3-pattern with no dashes and whose
$k$ leftmost or $k$ rightmost letters form either an increasing or
decreasing subword. Moreover, we find a bijection between
$n$-permutations that avoid  the pattern $132$ and begin with the
pattern $12$ and increasing rooted trimmed trees with $n+1$ nodes.
\end{abstract}

\maketitle
\thispagestyle{empty}

\section{Introduction and Background}

All permutations in this paper are written as words $\pi=a_1 a_2\cdots
a_n$, where the $a_i$ consist of all the integers $1,2,\ldots,n$.

A {\em pattern} is a word on some alphabet of letters, where some of the
letters may be separated by dashes. In our notation, the classical
permutation
patterns, first studied systematically by Simion and Schmidt
\cite{SchSim}, are of the form $p~=~1~-~3~-~2$, the dashes indicating
that
the letters in a permutation corresponding to an occurrence of $p$
do not have to be adjacent. In the classical case, an occurrence of a
pattern $p$ in a permutation $\pi$ is a subsequence in $\pi$ (of the
same length as the length of $p$) whose letters are in the same
relative order as those in~$p$.  For example, the permutation~264153
has only one occurrence of the pattern $1-2-3$, namely the
subsequence~245. Note that a classical pattern should, in our notation,
have dashes at
the beginning and end. Since most of the patterns considered in this
paper
satisfy this, we suppress these dashes from the notation.

In \cite{BabStein} Babson and Steingr\'{\i}msson introduced {\em
  generalized permutation patterns (GPs)} where two adjacent letters
in a pattern may be required to be adjacent in the permutation. Such
an adjacency requirement is indicated by the absence of a dash between
the corresponding letters in the pattern. Thus, a
pattern with no dashes corresponds to a contiguous subword anywhere in
a permutation. For example, the permutation $\pi=516423$ has only one
occurrence of the pattern $2$-$31$, namely the subword 564, but the
pattern $2$-$3$-$1$ occurs
also in the subwords 562 and 563. The motivation for introducing these
patterns in \cite{BabStein} was the study of Mahonian statistics.

A number of interesting results on GPs were obtained by Claesson
\cite{Claes}. Relations to several well studied combinatorial
structures, such as set partitions, Dyck paths, Motzkin paths and
involutions, were shown there. In \cite{Kit1} the present author
investigated simultaneous avoidance of two or more 3-letter GPs with
no dashes. Also there is a number of works concerning GPs by Mansour
(see for example \cite{Mans1, Mans2}).

In this paper we consider avoidance some generalized 3-patterns with
additional restrictions. The restrictions consist of demanding that a
permutation begin or end with the pattern $12\ldots k$ or the pattern
$k(k-1)\ldots 1$.

It turns out that the number of permutations that avoid the pattern
$1-32$ and end with the pattern $12 \ldots k$ is a linear combination of
the Bell numbers. The $n$-th Bell number is the number of ways a set of
$n$ elements can be partitioned into nonempty subsets. We find a
bijection between these permutations and all partitions of an
$(n-1)$-element set with one subset marked that satisfy certain special
conditions. In particular, in Theorem~\ref{correspondence1}, we
investigate the case $k=2$. We get that the total number of partitions
of an $(n-1)$-element set with one part marked, is equal to the number
of $(1-32)$-avoiding $n$-permutations that end with a $12$-pattern.
Lemma~\ref{theBinomialIdentity} gives us an identity involving the Bell
numbers and the Stirling numbers of the second kind, which seems to be
new. In Theorem~\ref{trees} we prove that the number of $132$-avoiding
$n$-permutations that begin with the pattern $12$ is equal to the number
of increasing rooted trimmed trees with $n+1$ nodes.

In Sections 4 -- 7, we give a complete description (in terms of {\em
exponential generating functions (e.g.f.)}) for the number of
permutations that avoid a pattern of the form $xyz$ and begin or end
with the pattern $12\ldots k$ or the pattern $k(k-1)\ldots 1$. We record
all the results concerning these e.g.f. in the table in Section~7. The
case $k=1$ is equivalent to the absence of the additional restriction.
This case was considered in \cite{ElizNoy} and \cite{Kit2}.

We observe that avoidance of some pattern with the additional
restrictions described above, in fact is equivalent to simultaneous
avoidance of several patterns. For example, beginning with the pattern
$12$ is equivalent to the avoidance of the pattern $[21)$ in the 
Babson-Steingr\'{\i}msson notation. Thus avoidance of the pattern
$132$ and beginning with the pattern $12$ is equivalent to
simultaneous avoidance of the patterns $132$ and $[21)$. Also, ending
with the pattern $123$ is equivalent to simultaneously avoiding the
patterns $(132]$, $(213]$, $(231]$, $(312]$ and $(321]$.

\section{Set partitions and pattern avoidance}

We recall some basic definitions.

A {\em partition} of a set $S$ is a family, $\pi = \{ A_1, A_2, \ldots ,
A_k \}$, of pairwise disjoint non-empty subsets of $S$ such that $S =
\cup_{i} A_i$. The total number of partitions of an $n$-element set is
called a {\em Bell number} and is denoted~$B_n$.

The {\em Stirling number of the second kind} $S(n,k)$ is the number of
ways a set with $n$ elements can be partitioned into $k$ disjoint,
non-empty subsets.

\begin{prop}\label{formulaForP(n,k)} Let $P(n,k)$ be the number of
$n$-permutations that avoid the pattern $1-32$ and end with the pattern
$12\ldots k$. Then $$P(n,k)=\displaystyle\sum_{i=0}^{n-k}{n-1 \choose
i}B_i.$$ \end{prop}

\begin{proof} Suppose a permutation $\pi=\sigma 1 \tau$ avoids the
pattern $1-32$ and ends with the pattern $12\ldots k$. The letters of
$\tau$ must be in  increasing order, since otherwise we have an
occurrence of the pattern $1-32$ involving 1. Also, $\sigma$ must avoid
$1-32$. If $|\sigma|=i$ then obviously $0\leq i \leq n-k$ and we can
choose the letters of $\sigma$ in ${n-1 \choose i}$ ways. By
\cite[Proposition 5]{Claes}, the number of $i$-permutations that avoid
the pattern $1-32$ is equal to $B_i$, hence there are $B_i$ ways to form
$\sigma$. \end{proof}

\begin{lemma}\label{theBinomialIdentity} We have \
$\displaystyle\sum_{i=0}^{n-1}{n \choose
i}B_i=\displaystyle\sum_{i=0}^{n}i\cdot S(n,i).$ \end{lemma}

\begin{proof} The identity can be proved from the recurrences for
$S(n,k)$ and $B_n$, but we give a combinatorial proof.

The left-hand side of the identity is the number of ways to choose $i$
elements from an $n$-element set, and then to make all possible
partitions of the chosen elements.

The right-hand side is the number of ways to partition a set with $n$
elements into $i$ disjoint non-empty subsets ($i=1,2, \ldots ,n$) and
mark  one of the subsets. For example if $n=4$ then $\overline{1}- 24 -
3$ and $1-\overline{24}-3$ are two different partitions, where the
marked subset is overlined.

A bijective correspondence between these combinatorial interpretations
is given by the following: For the left-hand side, after partitioning
the $i$ chosen elements, let the remaining $n-i$ elements form the
marked subset in the partition. \end{proof}

The formula for $P(n,k)$ in Proposition~\ref{formulaForP(n,k)},
applied to $k = 2$, and Lemma~\ref{theBinomialIdentity} now give the
following theorem:

\begin{thm}\label{correspondence1} The total number of partitions of an
$(n-1)$-element set with one part marked, is equal to the number of
$(1-32)$-avoiding $n$-permutations that end with the pattern $12$.
\end{thm}

We give now a direct combinatorial proof of this theorem.

\begin{proof} Suppose $P=S_1-S_2- \cdots -S_k$ is a partition of an
$(n-1)$-element set into $k$ subsets with one marked subset and $T_i$ is
the word that consists of all elements of $S_i$ in increasing order. We
may, without loss of generality, assume that $\min(S_i)<\min(S_j)$ if 
$i>j$. In particular, $1 \in S_k$. There are two cases possible:

\begin{itemize}
\item[1)] $S_k=\{1\}$ ($S_k$ is not marked set);
\item[2)] Either $S_k=\overline{1}$ or $1 \in S_k$ and $|S_k|\geq 2$.
\end{itemize}

In the first case, to a partition $P=S_1-S_2- \cdots
-\overline{S_i}-\cdots -S_{k-1}-1$ we associate the permutation
$\pi(P)=nT_1T_2 \ldots T_{i-1}T_{i+1} \ldots T_{k-1}1T_i$, which is
$(1-32)$-avoiding and ends with the pattern $12$ since $S_i \neq
\emptyset$. For example $4-\overline{23}-1 \mapsto 54123$.

In the second case we adjoin $n$ to a marked subset, and then consider
the permutation $\pi(P)=T_1T_2 \ldots T_k$. This permutation is
obviously $(1-32)$-avoiding since $\min(S_i)<\min(S_j)$ if $i>j$ and the
letters in $T_i$ are in increasing order. Also it ends with the pattern
$12$. For example $5-\overline{34}-12 \mapsto 534612$, and
$5-234-\overline{1} \mapsto 523416$.

Obviously in both cases we have an injection.

Now it is easy to see that the correspondence above is a surjection as
well. Indeed, for any $(1-32)$-avoiding permutation $\pi$ that ends with
the pattern $12$, we can check if $\pi$ begins with $n$ or not and
according to this we have either case 1) or 2). In the first case, we
remove $n$, then read $\pi$ from left to right and consider all maximal
increasing intervals. The elements of each such interval correspond to
some subset, and we let all the letters to the right of 1 constitute the
marked subset. In the second case, we divide $\pi$ into maximal
increasing intervals, and let the letters of each interval correspond to
a subset. Then we let the interval containing $n$ be the marked subset.
Thus we have a surjection.
So the correspondence is a bijection and the theorem is proved.
\end{proof}

The following theorem generalizes Theorem~\ref{correspondence1}.

\begin{thm}\label{correspondence2} Let $P=S_1-S_2- \cdots -S_{\ell}$ be
a partition of $\{ 1,2, \ldots, n-1 \}$ into $\ell$ subsets with subset
$S_i$ marked. We assume also that $1 \in S_{\ell}$. Then $P(n,k)$ counts
all possible marked partitions of $\{ 1,2, \ldots, n-1 \}$ that satisfy
the following conditions:
\begin{itemize}
\item[1)] if $i=\ell$ (the last subset is marked) then $|S_{\ell}| \geq
k-1$;
\item[2)] if $i \neq \ell$ and $|S_{\ell}| \neq 1$ then $|S_{\ell}| \geq
k$;
\item[3)] if $i \neq \ell$ and $|S_{\ell}|=1$ then $|S_i| \geq k-1$.
\end{itemize} \end{thm}

\begin{proof} A proof of this theorem is similar to the proof of
Theorem~\ref{correspondence1}. We assume that $\min(S_i)<\min(S_j)$ for
$i>j$ and consider three cases.

If a partition satisfies 1), that is $P=S_1-S_2- \cdots
-\overline{S_{\ell}}$ and $|S_{\ell}| \geq k-1$, then adjoining $n$ to
$S_{\ell}$ guarantees that the permutation $\pi(P)=T_1T_2 \ldots
T_{\ell}$, which is $(1-32)$-avoiding, ends with $k$ letters in
increasing order.

In case 2), we adjoin $n$ to the marked subset and consider
$\pi(P)=T_1T_2 \ldots T_{\ell}$. This permutation avoids the pattern
$1-32$ and ends with the pattern $12\ldots k$ since $|S_{\ell}| \geq k$.

In case 3), to a partition $P=S_1-S_2- \cdots -\overline{S_i}-\cdots
-S_{k-1}-1$ we associate the permutation $\pi(P)=nT_1T_2 \ldots
T_{i-1}T_{i+1} \ldots T_{k-1}1T_i,$ which is $(1-32)$-avoiding and ends
with at least $k$ letters in increasing order since $|S_i| \geq k-1$.

That this correspondence is a bijection can be shown in a way similar to
the proof of Theorem~\ref{correspondence1}.  \end{proof}

\section{Increasing rooted trimmed trees and pattern avoidance}

In an {\em increasing rooted tree}, nodes are numbered and the numbers
increase as we move away from the root. A {\em trimmed tree} is a tree
where no node has a single leaf as a child (every leaf has a sibling).

\begin{thm}\label{trees} Let $A_n$ denote the set of all
$n$-permutations that avoid the pattern $132$ and begin with the pattern
$12$. The number of permutations in $A_n$ is equal to the number of
increasing rooted trimmed trees (IRTTs) with $n+1$ nodes. \end{thm}

\begin{proof} A {\em right-to-left minimum} of a permutation $\pi$ is
an element $a_i$ such that $a_i<a_j$ for every $j>i$.

We describe a bijective correspondence $F$ between the permutations in
$A_n$ and IRTTs with $n+1$ nodes.

Suppose $\pi\in A_n$ and $\pi=P_0a_0P_1a_1\ldots P_ka_k$, where $a_i$
are the right-to-left minima of $\pi$ and $P_j$ are (possibly empty) 
subwords of $\pi$. We construct a IRTT $T = F(\pi)$ with $n+1$ nodes as
follows. The root of $T$ is labelled by $0$ and $a_0, a_1,\ldots , a_k$
are the labels of the root's children if we read them from left to
right. Then we let the right-to-left minima of $P_i$ be the labels of
the children of $a_i$ and so on. It is easy to see that, since $\pi$ 
avoids $132$ and begins with $12$, $T$ avoids limbs of length~2. Also,
$T$ is an increasing rooted tree and hence $T$ is a IRTT. For instance,
$$F(2\; 9\; 10\; 5\; 3\; 1\; 11\; 13\; 14\; 8\; 12\; 7\; 4\; 6) =
  
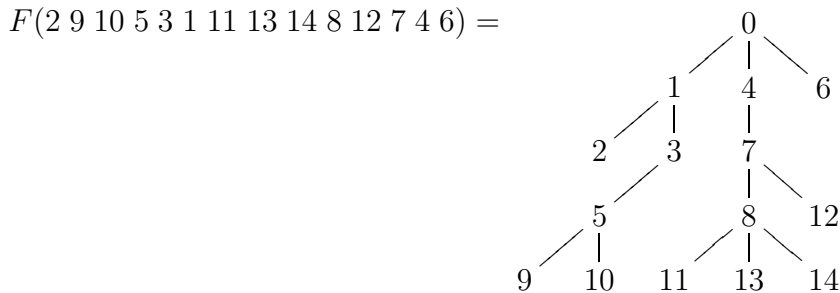
\xymatrix@!R@!C@=2ex@M=0.6ex{
      &        &         & 0\LDR   \\
      &        &  1\LD   & 4\DOWN     & 6   \\
      &   2    &  3\LEFT & 7\DR    \\
      &  5\LD  &         & 8\LDR   &   12 \\
   9  &  10    &  11     & 13      &   14 \\
}$$

Obviously, the correspondence $F$ is an injection.

To see, that $F$ is a surjection, we show how to construct the
permutation $\pi\in A_n$ that corresponds to a given IRTT $T$. The main
rule is the following: If $a_i$ and $a_j$ are siblings, and $a_i<a_j$,
then the labels of the nodes of the subtree below $a_j$, are
all the letters in $\pi$ between $a_i$ and $a_j$, that is, $a_{i+1}$,
$a_{i+2}, \ldots, a_{j-1}$. If $a_i$ is a single child, then the labels
of the nodes of the subtree below $a_i$ appear immediately left of $a_i$
in $\pi$. That is, if there are $k$ nodes in the subtree below $a_i$
then the $k$ corresponding labels form the subword 
$a_{i-k}a_{i-k+1}\ldots
a_{i-1}$. We now start from the first level of $T$, which consists of
the root's children, and apply this rule. After that we consider the
second level and so on. The fact that $T$ is a IRTT ensures that $\pi$
avoids the pattern $132$ and begins with the pattern~$12$. Thus, $F$ is
a bijection.
\end{proof}

\section{Avoiding $132$ and beginning with $12 \ldots k$ or
$k(k-1)\ldots 1$}

Let $E_{q}^{p}(x)$ denote the e.g.f. for the number of permutations that
avoid the pattern $q$ and begin with the pattern $p$.

If $k=1$, then there is no additional restriction, that is, we are
dealing with avoidance of the pattern $132$ (no dashes) and thus
\begin{equation}
E_{132}^{1}(x) = \frac{1}{1 - \int_{0}^{x}e^{-t^2/2}\ dt},
\label{formulaF1}
\end{equation}
since this result is a special case of \cite[Theorem 4.1]{ElizNoy} and
\cite[Theorem 12]{Kit2}.

\begin{thm}\label{132increasingBegin} We have
$$E_{132}^{12}(x) = \frac{e^{-x^2/2}}{1 - \int_{0}^{x}e^{-t^2/2}\ dt} -
x - 1,$$
and for $k \ge 3$
$$E^{12\ldots k}_{132}(x) =
E^1_{132}(x)\displaystyle\int_{0}^{x}\int_{0}^{t_{k-2}}\cdots
\int_{0}^{t_2}\left( e^{-t_1^2/2} - \frac{t_1+1}{E^1_{132}(t_1)} \right)
dt_1dt_2 \cdots dt_{k-2}.$$
\end{thm}

\begin{proof} Let $E_{n, k}$ denote the number of $n$-permutations that
avoid the pattern $132$ and begin with an increasing subword of length
$k > 0$. Let $\pi$ be such a permutation of length $n+1$. Also, suppose
$k \neq 2$. If $\pi = \sigma 1 \tau$ then either $\sigma = \epsilon$ or
$\sigma \neq \epsilon$ where $\epsilon$ denotes the empty word. If
$\sigma = \epsilon$ then $\tau$ must avoid $132$ and begin with an
increasing subword of length $k-1$. Otherwise $\sigma$ must avoid $132$
and begin with an increasing subword of length $k$, whereas $\tau$ must
begin with the pattern $12$, or be a single letter (there are $n$ ways
to choose this letter), or be $\epsilon$. This leads to the following:

\begin{equation}
E_{n+1,k} = E_{n,k-1} + \displaystyle\sum_{i\geq 0}{n \choose
i}E_{i,k}E_{n-i,2} + nE_{n-1,k} + E_{n,k}.
\label{equality1}
\end{equation}

Multiplying both sides of the equality with $x^n/n!$ and summing over
all $n$ we get the following differential equation
\begin{equation}
\frac{d}{d x}E^{12\ldots k}_{132}(x) = (E^{12}_{132}(x) + x +
1)E^{12\ldots k}_{132}(x) + E^{12\ldots (k-1)}_{132}(x),
\label{equality2}
\end{equation}
with the initial conditions $E^{12\ldots k}_{132}(0) = 0$ for $k\geq 3$.

Observe that equality~(\ref{equality2}) is not valid for $k=2$. Indeed,
if $k=2$, then it is incorrect to add the term $E_{n,k-1} = E_{n,1}$
in~(\ref{equality1}), since this term counts the number of permutations
$\pi = 1 \tau$ with the only restriction for $\tau$ that it must avoid
$132$. The absence of an additional restriction for $\tau$ means that
the 3 leftmost letters of $\pi$ could form the pattern $132$. However,
we can use~(\ref{equality2}) to find $E^{12}_{132}(x)$ by letting $k$
equal~1. In this case we have
$$\frac{d}{d x}E^1_{132}(x) = (E^{12}_{132}(x) + x + 1)E^1_{132}(x),$$
which gives
\begin{equation}
E^{12}_{132}(x) = \frac{e^{-x^2/2}}{1 - \int_{0}^{x}e^{-t^2/2}\ dt} - x
- 1.
\label{formulaF2}
\end{equation}

For the case $k\geq 3$, it is convenient to write $E^{12}_{132}(x)$ in
the form
$$E^{12}_{132}(x) = B^{\prime}(x)-x-1,$$
where $B(x)=-\ln (1-\int_{0}^{x}e^{-t^2/2}\ dt)$ and thus
$B^{\prime}(x)=\exp(B(x)-\frac{x^2}{2})$. So~(\ref{equality2}) is
equivalent to the differential equation
$$\frac{d}{d x}E^{12\ldots k}_{132}(x) = B^{\prime}(x)E^{12\ldots
k}_{132}(x) + E^{12\ldots (k-1)}_{132}(x)$$
which has the solution 
\bigskip\bigskip
\begin{eqnarray*}
&&E^{12\ldots k}_{132}(x) =
e^{B(x)}\displaystyle\int_{0}^{x}e^{-B(t)}E^{12\ldots
  (k-1)}_{132}(t)\;dt = \\[20pt]
&&E^1_{132}(x)\displaystyle\int_{0}^{x}\frac{E^{12\ldots
(k-1)}_{132}(t)}{E^{1}_{132}(t)}dt =\\[20pt]
&&E^1_{132}(x)\displaystyle\int_{0}^{x}\int_{0}^{t_2}\frac{E^{12\ldots
(k-2)}_{132}(t_1)}{E^{1}_{132}(t_1)}dt_1dt_2 =\\[20pt]
&& E^1_{132}(x)\displaystyle\int_{0}^{x}\int_{0}^{t_{k-2}}\cdots
\int_{0}^{t_2}\frac{E^{12}_{132}(t_1)}{E^1_{132}(t_1)}dt_1dt_2 \cdots
dt_{k-2}.
\end{eqnarray*}

\bigskip
\noindent
Using (\ref{formulaF1}) and (\ref{formulaF2}) we now get the desired
result.
\end{proof}

Using the formula for $E^{12\ldots k}_{132}(x)$ in
Theorem~\ref{132increasingBegin} one can derive, in particular, that
$$E^{123}_{132}(x) = -\frac{1}{2} - x - \frac{x^2}{2} + \frac{\left(
1+\frac{x}{2}\right)e^{-x^2/2} - \frac{1}{2}}{1-\int_{0}^{x}e^{-t^2/2}\
dt}.$$

\begin{thm}\label{132decreasingBegin} For $k \geq 2$
$$E_{k(k-1)\ldots 1}^{132}(x) =
\frac{E_1^{132}(x)}{(k-1)!}\displaystyle\int_{0}^{x}t^{k-1}e^{-t^2/2}\
dt.$$ \end{thm}

\begin{proof} We proceed as in the proof of
Theorem~\ref{132increasingBegin}.

Let $R_{n, k}$ denote the number of $n$-permutations that avoid the
pattern $132$ and begin with a decreasing subword of length $k > 1$ and
let $\pi$ be such a permutation of length $n+1$. Suppose also that $\pi
= \sigma 1 \tau$. If $\tau = \epsilon$ then, obviously, there are
$R_{n,k}$ ways to choose $\sigma$. If $|\tau| = 1$, that is, 1 is in the
second position from the right in $\pi$, then there are $n$ ways to
choose the rightmost letter in $\pi$ and we multiply this by
$R_{k,n-1}$, which is the number of ways to choose $\sigma$. If $|\tau|
> 1$ then $\tau$ must begin with the pattern $12$, otherwise the letter
1 and the two leftmost letters of $\tau$ form the pattern $132$, which
is forbidden. So, in this case there are $\sum_{i \geq 0}{n \choose i}R_{i,k}E_{n-i,2}$ such permutations with
the right properties, where $i$ indicates the length of $\sigma$ and
$E_{n-i,2}$ is defined in the proof of
Theorem~\ref{132increasingBegin}. In the last formula, of course,
$R_{i,k} = 0$ if $i<k$. Finally we have to consider the situation when
1 is in the $k$-th position. In this case we can choose the letters of
$\sigma$ in ${n \choose k-1}$ ways, write them in decreasing order and
then choose $\tau$ in $E_{n-k+1,2}$ ways. Thus

\begin{equation}
R_{n+1,k} = R_{n,k} + nR_{n-1,k} + \displaystyle\sum_{i\geq 0}{n \choose
i}R_{i,k}E_{n-i,2} + {n \choose k-1}E_{n-k+1,2}.
\label{equality3}
\end{equation}

We observe that (\ref{equality3}) is not valid for $n=k-1$ and $n=k$.
Indeed, if 1 is in the $k$-th position in these cases, the term ${n
\choose k-1}E_{n-k+1,2}$, which counts the number of such permutations,
is zero, whereas there is one ``good'' $(n+1)$-permutation in the case
$n=k-1$ and $n$ ``good'' $(n+1)$-permutations in case $n=k$. Multiplying
both sides of the equality with $x^n/n!$, summing over $n$ and using the
observation above (which gives the term $x^{k-1}/(k-1)! + kx^k/k!$ in
the right-hand side of Equalion~(\ref{equality4})), we get

\begin{equation}
\frac{d}{d x} E^{k(k-1)\ldots 1}_{132}(x) = (E^{12}_{132}(x) + x +
1)\left(E^{k(k-1)\ldots 1}_{132}(x)+ \frac{x^{k-1}}{(k-1)!}\right),
\label{equality4}
\end{equation}
with the initial condition $E_{k(k-1)\ldots 1}^{132}(0) = 0$. We solve
the equation in the way proposed in Theorem~\ref{132increasingBegin} and
get
$$E^{k(k-1)\ldots 1}_{132}(x) =
\frac{E^1_{132}(x)}{(k-1)!}\int_{0}^{x}\frac{(E^{12}_{132}(t) + t +
1)t^{k-1}}{E^1_{132}(t)}\ dt =
\frac{E^1_{132}(x)}{(k-1)!}\int_{0}^{x}t^{k-1}e^{-t^2/2}\ dt.$$
\end{proof}

For instance,
$$E^{21}_{132}(x) = \frac{1-e^{-x^2/2}}{1-\int_{0}^{x}e^{-t^2/2}\ dt}
\mbox{\ \ \ and\ \ \ } E^{321}_{132}(x) = \frac{1}{2}\left(-1 +
\frac{1-xe^{-x^2/2}}{1-\int_{0}^{x}e^{-t^2/2}\ dt}\right).$$

Moreover, the integral $\int_{0}^{x}t^{k-1}e^{-t^2/2}\ dt$ from the
formula for $E^{k(k-1)\ldots 1}_{132}(x)$ can be solved to show that
$E^{k(k-1)\ldots 1}_{132}(x)$ equals 
$$
\frac{(k/2-1)!2^{k/2-1}}{(k-1)!(1-\sqrt{\frac{\pi}{2}}\erf(x))}\left( 1
- e^{-x^2/2} \displaystyle\sum_{i=0}^{k/2-1}\frac{x^{2i}}{2^ii!}\right),
$$
if $k$ is even, and 
$$
\frac{1}{(k-1)!!}\left( -1 + \frac{1}{1-\sqrt{\frac{\pi}{2}}\erf(x)}
\left( 1 - e^{-x^2/2}
\displaystyle\sum_{i=0}^{(k-3)/2}\frac{x^{2i+1}}{(2i+1)!!}
\right)\right)
$$
if $k$ is odd.

In the formula above, $\erf(x)$ is the {\em error function}:$$\erf(x) =
\frac{2}{\sqrt{\pi}}\int_{0}^{x}e^{-t^2}\ dt.$$

\section{Avoiding $123$ and beginning with $k(k-1) \ldots 1$ or $12
\ldots k$}

If $k=1$, we have no additional restrictions and, according to
\cite[Theorem 4.1]{ElizNoy},
$$E^1_{123}(x) =
\frac{\sqrt{3}}{2}\frac{e^{x/2}}{\cos\left(\frac{\sqrt{3}}{2}x +
\frac{\pi}{6}\right)}.$$

\begin{thm}\label{123decreasingBegin} For $k \geq 2$
$$E^{k(k-1)\ldots 1}_{123}(x) =
\frac{e^{x/2}}{(k-1)!\cos\left(\frac{\sqrt{3}}{2}x +
\frac{\pi}{6}\right)}\int_{0}^{x}e^{-t/2}t^{k-1}\sin\left(\frac{\sqrt{3}}{2}t
+ \frac{\pi}{6}\right)\ dt.$$

In particular, $$E^{21}_{123}(x) =
\frac{\sqrt{3}}{2}\tan\left(\frac{\sqrt{3}}{2}x + \frac{\pi}{6}\right)
-x -\frac{1}{2}.$$ \end{thm}

\begin{proof} Let $P_{n, k}$ denote the number of $n$-permutations that
avoid the pattern $123$ and begin with a decreasing subword of length
$k$. We observe that we can use arguments similar to the proof of
Theorem~\ref{132decreasingBegin} to get the recurrence formula for
$P_{n,k}$. Indeed, we only need to write the letter $P$ instead of $R$
and $E$ in~(\ref{equality3}):

\begin{equation}
P_{n+1,k} = P_{n,k} + nP_{n-1,k} + \displaystyle\sum_{i\geq 0}{n \choose
i}P_{i,k}P_{n-i,2} + {n \choose k-1}P_{n-k+1,2}.
\label{equality5}
\end{equation}

This formula is valid for $k>1$. Multiplying both sides of the equality
with $x^n/n!$, summing over $n$ and reasoning as in the proof of
Theorem~\ref{132decreasingBegin}, we get:

\begin{equation}
\frac{d}{d x} E^{k(k-1)\ldots 1}_{123}(x)= (E^{21}_{123}(x) + x +
1)\left(E^{k(k-1)\ldots 1}_{123}(x) + \frac{x^{k-1}}{(k-1)!}\right),
\label{equality5.1}
\end{equation}
with the initial condition $E^{k(k-1)\ldots 1}_{123}(0) = 0$. To
solve~(\ref{equality5.1}), we need to know $E^{21}_{123}(x)$. To find
it, we consider the case $k=1$. In this case we have almost the same
recurrence as we have in~(\ref{equality5}), but we must remove the last
term in the right-hand side:
$$P_{n+1,1} = P_{n,1} + nP_{n-1,1} + \displaystyle\sum_{i\geq 0}{n
\choose i}P_{i,k}P_{n-i,2}.$$

After multiplying both sides of the last equality with $x^n/n!$ and
summing over $n$, we have
$$\frac{d}{d x} E^1_{123}(x) = (E^{21}_{123}(x) + x + 1)E^1_{123}(x),$$
and thus
$$E^{21}_{123}(x)=\frac{\frac{d}{d x} E^1_{123}(x)}{P_1(x)}-x-1=
\frac{\sqrt{3}}{2}\tan\left(\frac{\sqrt{3}}{2}x + \frac{\pi}{6}\right)
-x - \frac{1}{2}.$$

Now we solve~(\ref{equality5.1}) in the way we solved
Equalion~(\ref{equality4}) and get
$$E^{k(k-1)\ldots 1}_{123}(x) =
\frac{e^{x/2}}{(k-1)!\cos\left(\frac{\sqrt{3}}{2}x +
\frac{\pi}{6}\right)}\int_{0}^{x}e^{-t/2}t^{k-1}\sin\left(\frac{\sqrt{3}}{2}t
+ \frac{\pi}{6}\right)\ dt.$$
\end{proof}

The following theorem is straightforward to prove.

\begin{thm}\label{123increasingBegin} We have $E^{12\ldots
k}_{123}(x)=0$ for $k \geq 3$ and
$$E^{12}_{123}(x) = E^1_{123}(x) - E^{21}_{123}(x) =
\frac{\sqrt{3}}{2}\frac{e^{x/2}}{\cos\left(\frac{\sqrt{3}}{2}x +
\frac{\pi}{6}\right)} - \frac{1}{2} -
\frac{\sqrt{3}}{2}\tan\left(\frac{\sqrt{3}}{2}x +
\frac{\pi}{6}\right).$$ \end{thm}

\section{Avoiding $213$ and beginning with $k(k-1) \ldots 1$ or $12
\ldots k$}

If $k=1$, then by \cite[Theorem 4.1]{ElizNoy} or \cite[Theorem 12]{Kit2}

$$E^1_{312}(x) = \frac{1}{1 - \int_{0}^{x}e^{-t^2/2}\ dt}.$$

\begin{thm}\label{213increasingBegin} For $k\geq 2$
$$E^{12\ldots
k}_{213}(x)=\int_{0}^{x}\int_{0}^{t}\frac{s^{k-2}e^{T(t)-T(s)}}{(k-2)!(1-\int_{0}^{t}e^{-m^2/2}d
m)}\ dsdt,$$ where
$T(x)=-x^2/2+\int_{0}^{x}\frac{e^{-t^2/2}}{1-\int_{0}^{t}e^{-s^2/2}d s}\
d t$. \end{thm}

\begin{proof} Let $A_n$ denote the number of $n$-permutations that avoid
the pattern $213$ and let $B_n$ denote the number of $n$-permutations
that avoid $213$ and begin with the pattern $12\ldots k$. Let $C_n$
denote the number of $n$-permutation that avoid $213$, begin with the
pattern $12\ldots k$ and end with the pattern $12$ and let $D_n$ denote
the number of $n$-permutations that avoid $213$ and end with the pattern
$12$. Also, let $A(x)$, $B(x)$, $C(x)$ and $D(x)$ denote the e.g.f. for
the numbers $A_n$, $B_n$, $C_n$ and $D_n$ respectively.

We observe, that
$$D(x)=E^{12}_{132}(x)=e^{-x^2/2}/(1-\int_{0}^{x}e^{-t^2/2}d t) - x -
1,$$ since, by using the reverse and complement discussed in the next
section, there are as many permutations that avoid the pattern $213$ and
end with the pattern $12$ as those that avoid the pattern $132$ and
begin with the pattern $12$. Also, $A(x)=E^1_{312}(x)$ and $B(x)=
E^{12\ldots k}_{213}(x)$.

Suppose now that $\pi = \sigma (n+1) \tau$ is an $(n+1)$-permutation
that avoids the pattern $213$ and begins with the pattern $12\ldots k$.
So $\sigma$ must avoid $213$, begin with $12\ldots k$, but also end with
the pattern $12$ since otherwise the two rightmost letters of $\sigma$
together with the letter $(n+1)$ form the pattern $213$, which is
forbidden. For $\tau$, there is only one restriction --- avoidance of 
$213$. So if $|\sigma|=i$ then we can choose the letters of $\sigma$ in
${n \choose i}$ ways, which gives $\sum_{i\geq 0}{n \choose
i}C_iA_i$ permutations that avoid the pattern $213$ and begin with the
pattern $12\ldots k$. Moreover, it is possible for $(n+1)$ to be in the
$k$th position, in which case we choose the letters of $\sigma$ in ${n
\choose k-1}$ ways and arrange them in increasing order. Thus
$$B_{n+1} = \displaystyle\sum_{i\geq 0}{n \choose i}C_iA_{n-i} + {n
\choose k-1}A_{n-(k-1)}.$$

Multiplying both sides of this equality with $x^n/n!$ and summing 
over~$n$, we get
\begin{equation}
B^{\prime}(x) = \left(C(x)+\frac{x^{k-1}}{(k-1)!}\right)A(x),
\label{equality6}
\end{equation}
with the initial condition $B(0)=0$.

To solve (\ref{equality6}) we need to find $C(x)$. Let $\pi = \sigma
(n+1) \tau$ be an $(n+1)$-permutation that avoids the pattern $213$,
begins with the pattern~$12\ldots k$ and ends with the pattern $12$.
Reasoning as above, $\sigma$ must avoid the pattern $213$, begin with
the pattern $12\ldots k$ and end with the pattern $12$, whereas $\tau$ 
must avoid $213$ and end with the pattern~$12$. This gives
$\sum_{i\geq 0}{n \choose i}C_iD_{n-i}$ permutations
counted by $C_{n+1}$. Also, the letter $(n+1)$ can be in the $k$th
position, which gives ${n \choose k-1}D_{n-(k-1)}$ permutations, and
this letter can be in the $(n+1)$st position, which gives $C_n$
permutations that avoid the pattern $213$, begin with the pattern
$12\ldots k$ and end with the pattern $12$. Also, if $n+1=k$ and all the
letters are arranged in increasing order, then $(n+1)$ is in the
$(n+1)$st position, but this permutation is not counted by $C_n$ above.
So
$$C_{n+1} = \displaystyle\sum_{i\geq 0}{n \choose i}C_iD_{n-i} + {n
\choose k-1}D_{n-(k-1)} + C_n + {\delta}_{n,k-1},$$
where $\delta_{n,k}$ is the Kronecker delta, that is,

\[ \delta_{n,k} = \left\{ \begin{array}{ll} 1, & \mbox{if $n=k$,}
\\ 0, & \mbox{else.}
\end{array}
\right. \]

Multiplying both sides of the equality with $x^n/n!$ and summing
over~$n$, we get

\begin{equation}
C^{\prime}(x) = (D(x)+1)C(x) + (D(x)+1)\frac{x^{k-1}}{(k-1)!}.
\label{equality7}
\end{equation}

To solve (\ref{equality7}), it is convenient to introduce the function
$T(x)$ such that $T^{\prime}(x)=D(x)+1$. Thus
$$T(x) = x + \int_{0}^{x}D(t)dt =
-x^2/2+\int_{0}^{x}\frac{e^{-t^2/2}}{1-\int_{0}^{t}e^{-s^2/2}d s}\ d
t,$$
and we need to solve the equation
$$C^{\prime}(x)= T^{\prime}(x)C(x) +
T^{\prime}(x)\frac{x^{k-1}}{(k-1)!},$$
with $C(0) = 0$.

The solution to this equation is given by
$$C(x)=e^{T(x)}\int_{0}^{x}e^{-T(t)}T^{\prime}(t)\frac{t^{k-1}}{(k-1)!}d
t =
-\frac{x^{k-1}}{(k-1)!}+e^{T(x)}\int_{0}^{x}e^{-T(t)}\frac{t^{k-2}}{(k-2)!}d
t.$$

Now we substitute $C(x)$ into~(\ref{equality6}) to get the desired
result. \end{proof}

\begin{thm} For $k\geq 2$
$$E^{k(k-1)\ldots 1}_{213}(x) = -\frac{x^{k-1}}{(k-1)!} +
\displaystyle\sum_{n=0}^{k-2}\int_{0}^{x}\int_{0}^{t_n}\cdots\int_{0}^{t_1}\frac{C_{k-n}(t)
+ \delta_{n,k-2}}{1-\int_{0}^{t}e^{-m^2/2}d m}d td t_1 \cdots dt_n,$$
where
$$C_k(x) =
e^{T(x)}\displaystyle\int_{0}^{x}\int_{0}^{t_{k-2}}\cdots\int_{0}^{t_1}
e^{-T(t)}\left( \frac{e^{-t^2/2}}{1-\int_{0}^{t}e^{-m^2/2}dm} - t -
1\right)dt dt_1 \cdots dt_{k-2},$$
with
$T(x)=-x^2/2+\int_{0}^{x}\frac{e^{-t^2/2}}{1-\int_{0}^{t}e^{-s^2/2}d s}\
d t$.
\end{thm}

\begin{proof} Let $A_n$ denote the number of $n$-permutations that avoid
the pattern $213$ and let $B_{n,k}$ denote the number of
$n$-permutations that avoid $213$ and begin with the pattern
$k(k-1)\ldots 1$ for $k\geq 2$. Let $C_{n,k}$ denote the number of
$n$-permutation that avoid $213$, begin with $k(k-1)\ldots 1$ for $k\geq
2$ and end with the pattern $12$ and let $D_n$ denote the number of
$n$-permutations that avoid $213$ and end with the pattern $12$. Also,
let $A(x)$, $B_k(x)$, $C_k(x)$ and $D(x)$ denote the e.g.f. for the
numbers $A_n$, $B_{n,k}$, $C_{n,k}$ and $D_n$ respectively. In the proof
of Theorem~\ref{213increasingBegin} it was shown that
$D(x)=e^{-x^2/2}/(1-\int_{0}^{x}e^{-t^2/2}d t) - x - 1$ and
$A(x)=E^1_{312}(x)$. Moreover, $B_k(x)= E^{k(k-1)\ldots 1}_{213}(x)$.

Suppose now that $\pi = \sigma (n+1) \tau$ is an $(n+1)$-permutation
that avoids $213$ and begins with the pattern $k(k-1)\ldots 1$. So
$\sigma$ must avoid $213$, begin with $k(k-1)\ldots 1$, but also end
with the pattern $12$ if $|\sigma|\geq 2$, since otherwise the two
rightmost letters of $\sigma$ together with the letter $(n+1)$ form the
pattern $213$ which is forbidden. For $\tau$, there is only one
restriction - avoidance of $213$. So if $|\sigma|=i$ then we can choose
the letters of $\sigma$ in ${n \choose i}$ ways, which gives
$\sum_{i\geq 0}{n \choose i}C_{i,k}A_i$ permutations
counted by $B_{n+1,k}$. Also, it is possible for $(n+1)$ to be the
leftmost letter, in which case the remaining letters must form a
$n$-permutation that avoids $213$ and begins with the pattern
$(k-1)(k-2)\ldots 1$. Thus
\begin{equation}
B_{n+1,k} = \displaystyle\sum_{i\geq 0}{n \choose i}C_{i,k}A_{n-i} +
B_{n,k-1}.
\label{equality9}
\end{equation}
However, this formula is not valid when $k=2$ and $n=0$. Indeed, since
$B_{0,1} = A_{0}=1$, it follows from the formula that $B_{1,2}=1$, which
is not true, since $B_{1,2}$ must be $0$. So, in the right-hand side
of~(\ref{equality9}), we need to subtract the term
\[ \gamma_{n,k} = \left\{ \begin{array}{ll} 1, & \mbox{if $n=0$ and
$k=2$,}
\\ 0, & \mbox{else.}
\end{array}
\right. \]
Multiplying both sides of the obtained equality by $x^n/n!$ and summing
over $n$, we get, that for $k\geq 3$
\begin{equation}
\frac{d}{d x} B_k(x) = C_k(x)A(x) + B_{k-1}(x),
\label{equality10}
\end{equation}
with the initial condition $B_k(0)=0$, and
\begin{equation}
\frac{d}{d x} B_2(x) = C_2(x)A(x) + B_{1}(x) - 1,
\label{equality11}
\end{equation}
with the initial condition $B_2(0)=0$.

The solution to differential equations~(\ref{equality10})
and~(\ref{equality11}) is given by
$$B_k(x) = -\frac{x^{k-1}}{(k-1)!} +
\displaystyle\sum_{n=0}^{k-2}\int_{0}^{x}\int_{0}^{t_n}\cdots\int_{0}^{t_1}\frac{C_{k-n}(t)
+ \delta_{n,k-2}}{1-\int_{0}^{t}e^{-m^2/2}d m}d td t_1 \cdots dt_n.$$

So, to prove the theorem, we only need to find $C_k(x)$.

Suppose $\pi = \sigma (n+1) \tau$ be an $(n+1)$-permutation that avoids
the pattern $213$, begins with the pattern $k(k-1)\ldots 1$ and ends
with the pattern $12$. It is clear that $\sigma$ must avoid $213$, begin
with the pattern $k(k-1)\ldots 1$ and end with the pattern $12$, whereas
$\tau$ must avoid $213$ and end with the pattern $12$. There are $\sum_{i\geq
0}{n \choose i}C_{i,k}D_{n-i}$ permutations with these properties. Also,
the letter $(n+1)$ can be in the leftmost position, which gives
$C_{n,k-1}$ permutations, and $(n+1)$  can be in the rightmost 
position, which gives $C_{n,k}$ permutations, since in this case, two
letters immediately to the left of $(n+1)$ cannot form a descent. So,
$$C_{n+1,k} = \displaystyle\sum_{i\geq 0}{n \choose i}C_{i,k}D_{n-i} +
C_{n,k-1} + C_{n,k}.$$
Multiplying both sides of the equality with $x^n/n!$ and summing over
$n$, we get the following differential equation
\begin{equation}
C^{\prime}_k(x) = (D(x)+1)C_k(x) + C_{k-1}(x).
\label{equality20}
\end{equation}

As when solving Equation~(\ref{equality7}), it is convenient to
introduce the function $T(x)$ such that $T^{\prime}(x)=D(x)+1$.
Moreover, Equation~(\ref{equality20}) is similar to
Equation~(\ref{equality2}) and we can solve it in the same way. Also we
observe that from the definitions, $C_1(t)=D(t)$, and thus
$$C_k(x) =
e^{T(x)}\displaystyle\int_{0}^{x}\int_{0}^{t_{k-2}}\cdots\int_{0}^{t_1}
e^{-T(t)}C_1(t)dt dt_1 \cdots dt_{k-2}=$$
$$e^{T(x)}\displaystyle\int_{0}^{x}\int_{0}^{t_{k-2}}\cdots\int_{0}^{t_1}
e^{-T(t)}\left( \frac{e^{-t^2/2}}{1-\int_{0}^{t}e^{-m^2/2}dm} - t -
1\right)dt dt_1 \cdots dt_{k-2}.$$
\end{proof}

\section{Summarizing the results from sections 4,5 and 6}

We recall that the {\em reverse} $R(\pi)$ of a permutation $\pi=a_1a_2
\ldots a_n$ is the permutation $a_na_{n-1} \ldots a_1$ and the {\em
complement} $C(\pi)$ is the permutation $b_1b_2 \ldots b_n$ where
$b_i=n+1-a_i$. Also, $R \circ C$ is the composition of $R$ and $C$. We
call these bijections of $\mathcal{S}_n$ to itself {\em trivial}. Let
$\phi$ be an arbitrary trivial bijection. It is easy to see that, for
example, there are as many permutations avoiding the pattern 132 as
those avoiding the pattern $\phi(132)$. Moreover if, for instance, a
permutation $\pi$ begins with a decreasing pattern of length $k$, then
depending on $\phi$, $\phi(\pi)$ either begins with an increasing
pattern, or ends with either a decreasing or increasing pattern of
length $k$. This allows us to apply Theorems 6 -- 11 to a number of
other cases. We summarize all the obtained results concerning avoidance
of a generalized 3-pattern with no dashes and beginning or ending with
either increasing or decreasing subword, in the table below.

\begin{center}
\begin{tabular}{|c|c|c|c|c|}
\hline
 & {\bf avoid} & {\bf begin} & {\bf end} & {\bf e.g.f.} \\
\hline
        & 123 & $12\ldots k$ & -- &  \\
{\bf 1} & 123 & --           & $12\ldots k$ &
$\frac{\sqrt{3}}{2}\frac{e^{x/2}}{\cos(\frac{\sqrt{3}}{2}x +
\frac{\pi}{6})}$, if $k=1$ \\[4mm]
        & 321 & $k\ldots 21$ & -- &
$\frac{\sqrt{3}}{2}\frac{e^{x/2}}{\cos(\frac{\sqrt{3}}{2}x +
\frac{\pi}{6})} - \frac{1}{2} -
\frac{\sqrt{3}}{2}\tan(\frac{\sqrt{3}}{2}x + \frac{\pi}{6})$, if $k=2$
\\[4mm]
        & 321 & --               & $k\ldots 21$ & 0, if $k\geq 3$
\\[3mm]
\hline
        & 123 & $k\ldots 21$ & -- & \\
{\bf 2} & 123 &  --        & $k\ldots 21$  &
$\frac{\sqrt{3}}{2}\frac{e^{x/2}}{\cos(\frac{\sqrt{3}}{2}x +
\frac{\pi}{6})}$, if $k=1$ \\
        & 321 & $12\ldots k$ & -- &  \\
        & 321 & -- &  $12\ldots k$ &
$\frac{e^{x/2}\int_{0}^{x}e^{-t/2}t^{k-1}\sin(\frac{\sqrt{3}}{2}t +
\frac{\pi}{6}))\ dt}{(k-1)!\cos(\frac{\sqrt{3}}{2}x + \frac{\pi}{6})}$,
if $k \geq 2$ \\[4mm]
\hline
        & 132 & $12\ldots k$ & -- &  \\[3mm]
{\bf 3} & 213 & --           & $12\ldots k$ & $(1 -
\int_{0}^{x}e^{-t^2/2}\ dt)^{-1}$, if $k=1$ \\[4mm]
        & 312 & $k\ldots 21$ & -- & $e^{-x^2/2}(1 -
\int_{0}^{x}e^{-t^2/2}\ dt)^{-1} - x - 1$, if $k=2$ \\[4mm]
        & 231 & -- & $k\ldots 21$ & $(1 - \int_{0}^{x}e^{-t^2/2}\
dt)^{-1}\int_{0}^{x}\int_{0}^{t_{k-2}}\cdots \int_{0}^{t_2}(e^{-t_1^2/2}
- $ \\[3mm]
        &     &    &                  & $(t_1+1)(1 -
\int_{0}^{t_1}e^{-t^2/2}\ dt)) dt_1dt_2 \cdots dt_{k-2}$, if $k\geq 3$
\\[4mm]
\hline
        & 132 & $k\ldots 21$ & -- &  \\
{\bf 4} & 213 & --           &  $k\ldots 21$ & $(1 -
\int_{0}^{x}e^{-t^2/2}\ dt)^{-1}$, if $k=1$ \\[4mm]
        & 312 & $12\ldots k$ & -- & $\frac{1}{(k-1)!(1 -
\int_{0}^{x}e^{-t^2/2}\ dt)}\int_{0}^{x}t^{k-1}e^{-t^2/2}\ dt$, if
$k\geq 2$ \\
        & 231 & -- & $12\ldots k$ &  \\
\hline
        & 213 & $12\ldots k$ & -- & $(1 - \int_{0}^{x}e^{-t^2/2}\
dt)^{-1}$, if $k=1$ \\
{\bf 5} & 132 & --           & $12\ldots k$ &  \\
        & 231 & $k\ldots 21$ & -- &
$\int_{0}^{x}\int_{0}^{t}\frac{s^{k-2}e^{T(t)-T(s)}}{(k-2)!(1-\int_{0}^{t}e^{-m^2/2}d
m)}\ dsdt$ , if $k\geq 2$, where \\
        & 312 & -- & $k\ldots 21$ &
$T(x)=-x^2/2+\int_{0}^{x}\frac{e^{-t^2/2}}{1-\int_{0}^{t}e^{-s^2/2}d s}\
d t$\\[4mm]
\hline
        & 213 & $k\ldots 21$ & -- &  $(1 - \int_{0}^{x}e^{-t^2/2}\
dt)^{-1}$, if $k=1$  \\[4mm]
{\bf 6} & 132 & --           &  $k\ldots 21$ & $-\frac{x^{k-1}}{(k-1)!}
+
\displaystyle\sum_{n=0}^{k-2}$$\int_{0}^{x}\int_{0}^{t_n}\cdots\int_{0}^{t_1}\frac{C_{k-n}(t)
+ \delta_{n,k-2}}{1-\int_{0}^{t}e^{-m^2/2}d m}d td t_1 \cdots dt_n$,
\\[4mm]
        & 231 & $12\ldots k$ & -- & if $k\geq 2$, where
$C_k(x)=e^{T(x)}\int_{0}^{x}\int_{0}^{t_{k-2}}\cdots \int_{0}^{t_1}
e^{-T(t)}\cdot$ \\[4mm]
        & 312 & -- & $12\ldots k$ & $\left(
\frac{e^{-t^2/2}}{1-\int_{0}^{t}e^{-m^2/2}dm} - t - 1\right)dt dt_1
\cdots dt_{k-2}$ and $T(x)$ as above \\[4mm]
\hline
\end{tabular}
\end{center}

\end{document}